%&amstex          
\input amstex\documentstyle{amsppt}  
\pagewidth{12.5cm}\pageheight{19cm}\magnification\magstep1
\topmatter
\title Total positivity in the space of maximal tori\endtitle
\author G. Lusztig\endauthor
\address{Department of Mathematics, M.I.T., Cambridge, MA 02139}\endaddress
\thanks{Supported by NSF grant DMS-2153741}\endthanks
\endtopmatter   
\document
\define\mat{\matrix}
\define\endmat{\endmatrix}

\define\si{\sim}
\define\wt{\widetilde}

\define\op{\oplus}
   
\define\part{\partial}
\define\emp{\emptyset}

\define\m{\mapsto}
\define\do{\dots}

\define\sub{\subset}    

\define\T{\times}
\define\ti{\tilde}
\define\nl{\newline}
\redefine\i{^{-1}}

\define\a{\alpha}
\redefine\b{\beta}
\redefine\c{\chi}

\redefine\d{\delta}

\define\io{\iota}

\define\p{\pi}

\define\z{\zeta}
\define\x{\xi}

\define\CC{\bold C}

\define\RR{\bold R}

\define\TT{\bold T}

\define\cb{\Cal B}

\define\ce{\Cal E}

\define\ct{\Cal T}

\define\fg{\frak g}

\head Introduction\endhead
\subhead 0.1\endsubhead
Let $G$ be a connected reductive algebraic
group over $\CC$ with a fixed pinning
$(B^+,B^-,x_i,y_i (i\in I))$ as in \cite{L94}. Here $B^+,B^-$ are
opposed Borel subgroups of $G$ with unipotent radicals $U^+,U^-$
and $x_i:\CC@>>>U^+,y_i:\CC@>>>U^-$ are certain imbeddings; $I$ is a finite
set. Let $T=B^+\cap B^-$.

In \cite{L94} it was shown that each of $U^+,U^-,T,G$ contains a certain
sub-semigroup $U^+_{>0},U^-_{>0},T_{>0},G_{>0}$ which can be regarded as a
form over $\RR_{>0}$ and is called the totally positive part of the
corresponding object. (It is a cell in each case.) This generalized the
known definition in type $A$.

Let $\TT$ be the set of maximal tori of $G$ (a complex algebraic variety).
In this paper we define a subset $\TT_{>0}$ of $\TT$
(a cell) which can be regarded as a form over $\RR_{>0}$
of $\TT$. (See \S1.) We say that $\TT_{>0}$ is the totally positive part of $\TT$.
The definition of $\TT$ is based on considering both the pinning above
and a new pinning (the inverse of the original one).

\subhead 0.2\endsubhead
Let $W$ be the Weyl group of $G$ and let
$\{s_i;i\in I\}$ be the set of simple reflections of $W$.
Let $\nu$ be the number of reflections of $W$ and let $I_*$ be the set
of all $i_*=(i_1,i_2,\do,i_\nu)\in I^\nu$ such that
$s_{i_1}s_{i_2}\do s_{i_\nu}\in W$ has length $\nu$.
Let $i_*\in I_*$. Following \cite{L94} we denote by $U^+_{>0}$
(resp. $U^-_{>0}$) the image of the
(injective) map $\RR_{>0}^\nu@>>>U^+$ (resp. $\RR_{>0}^\nu@>>>U^-$)
given by
$$(a_1,a_2,\do,a_\nu)\m x_{i_1}(a_1)x_{i_2}(a_2)\do x_{i_\nu}(a_\nu)$$
(resp.
$$(a_1,a_2,\do,a_\nu)\m y_{i_1}(a_1)y_{i_2}(a_2)\do y_{i_\nu}(a_\nu).)$$

In \cite{L94,\S2} it is shown that $U^+_{>0}$ (resp. $U^-_{>0}$) is
independent of the choice of $i_*$ and it is
a sub-semigroup of $U^+$ (resp. $U^-$) said to be the totally
positive part of  $U^+$ (resp. $U^-$).

Let $T_{>0}$ be the subgroup of $T$ generated by $\z(\RR_{>0})$
for various homomorphisms of algebraic groups $\z:\CC^*@>>>T$;
this is said to be the totally positive part of $T$.

In \cite{L94,\S2} it is shown that $U^+_{>0}T_{>0}U^-_{>0}=U^-_{>0}T_{>0}
U^+_{>0}$ and that this subset
of $G$ (denoted by $G_{>0}$) is a sub-semigroup of $G$ said to be the
totally positive part of $G$.

Now $U^+_{>0},U^-_{>0},T_{>0},G_{>0}$ are cells of dimension
$\nu,\nu,r,2\nu+r$ (here $r$ is the rank of $G$).

\subhead 0.3\endsubhead
In \cite{L94,\S8} it is shown that
$$\{uB^+u\i;u\in U^-_{>0}\}=\{vB^-v\i;v\in U^+_{>0}\};$$
this subset of $\cb$ (denoted by $\cb_{>0}$)
is said to be the totally positive part of $\cb$.
Moreover, $u\m uB^+u\i$ (resp. $v\m vB^-v\i$) is a bijection
$U^-_{>0}@>>>\cb_{>0}$ (resp. $U^+_{>0}@>>>\cb_{>0}$) so that
$\cb_{>0}$ is a cell of dimension $\nu$.

\subhead 0.4\endsubhead
In \S2 we define a map $\p:G_{>0}@>>>\TT_{>0}$ and in \S5
we conjecture that this map is surjective. (The examples of $GL_2$ and
$GL_3$ are considered in \S3,\S4.)

\subhead 0.5\endsubhead
I thank Xuhua He for useful comments.

\head 1. Total positivity on $\cb^2$ and $\TT$\endhead
\subhead 1.1\endsubhead
We define a new pinning $(B^+,B^-,x'_i,y'_i (i\in I))$ for $G$ by
$x'_i(a)=x_i(-a)=x_i(a)\i$, $y'_i(a)=y_i(-a)=y_i(a)\i$ for
$i\in I,a\in\RR$. Let $U_{<0}^+,U_{<0}^-,T_{<0},G_{<0}$ be
the sub-semigroups of $U^+,U^-,T,G$ defined like
$U_{>0}^+,U_{>0}^-,T_{>0},G_{>0}$ but using the new pinning instead of
the one in 0.1. We have $T_{<0}=T_{>0}$ and

$U_{<0}^+=\{u\in U^+;u\i\in U_{>0}^+\}$,

$U_{<0}^-=\{u\in U^-;u\i\in U_{>0}^-\}$.
\nl
(We use that if $(i_1,i_2,\do,i_\nu)\in I_*$, then
$(i_\nu,i_{\nu-1},\do,i_1)\in I_*$.)
It follows that $G_{<0}=\{g\in G;g\i\in G_{>0}\}$.

Let $\cb_{<0}$ be the subset of $\cb$ 
defined like $\cb_{>0}$ in 0.3 but using the new pinning instead of
the one in 0.1. Thus we have
$$\cb_{<0}=\{u\i B^+u;u\in U_{>0}^-\}=\{v\i B^-v;v\in U_{>0}^+\}.$$
As in 0.3, $u\m u\i B^+u$ (resp. $v\m v\i B^-v$) is a bijection
$U^-_{>0}@>>>\cb_{<0}$ (resp. $U^+_{>0}@>>>\cb_{<0}$) so that
$\cb_{<0}$ is a cell of dimension $\nu$. Hence we have a bijection
$U^-_{<0}@>>>U^+_{<0}$, denoted by $u\m\ti{u}$, such that
$uB^+u\i=\ti{u}B^-\ti{u}\i$ for all $u\in U^-_{<0}$.

We set $\cb^2_{>0}=\cb_{>0}\T\cb_{<0}$ (a subset of $\cb^2=\cb\T\cb$).
We say that this is the totally positive part of $\cb^2$.

\proclaim{Proposition 1.2} If $(B,B')\in\cb^2_{>0}$, then $B,B'$ are
opposed, that is $B\cap B'\in\TT$.
\endproclaim
We have $B=uB^+u\i,B'=v\i B^+v$ for well defined
$u\in U_{>0}^-,v\in U_{>0}^-$. We have
$$B\cap B'=u(B^+\cap zB^+z\i)u\i=u(B^+\cap\ti{z}B^-\ti{z}\i)u\i$$
where $z=u\i v\i\in U_{<0}^-$ so that $\ti{z}\in U_{<0}^+$ is
defined. Hence
$$\align&B\cap B'=u\ti{z}(\ti{z}\i B^+\ti{z}\cap B^-)\ti{z}\i u\i
=u\ti{z}(B^+\cap B^-)\ti{z}\i u\i=u\ti{z}T\ti{z}\i u\i\\&
=u\wt{(u\i v\i)}T(\wt{(u\i v\i)})\i u\i\endalign$$
since $\ti{z}\in B^+$. This proves the proposition.

\proclaim{Proposition 1.3}The map $\io:\cb^2_{>0}@>\si>>\TT$ given by
$(B,B')\m B\cap B'$ is injective. (This map is well defined by 1.2.)
\endproclaim
Assume that $(B_1,B'_1)\in\cb^2_{>0}$,
$(B_2,B'_2)\in\cb^2_{>0}$ satisfy $B_1\cap B'_1=B_2\cap B'_2$.
Then $B_1,B'_1,B_2,B'_2$ contain a common maximal torus:
$B_1\cap B'_1=B_2\cap B'_2$.
By 1.2, $B_1,B'_1$ are opposed and $B_2,B'_1$ are opposed.
Since $B_1,B'_1,B_2$ contain a common maximal torus, it follows that
$B_1=B_2$.
Again, by 1.2, $B_1,B'_1$ are opposed and $B_1,B'_2$ are opposed.
Since $B_1,B'_1,B'_2$ contain a common maximal torus, it follows that
$B'_1=B'_2$. This proves the injectivity of $\io$.

\subhead 1.4\endsubhead
The image of the injective map $\io:\cb^2_{>0}@>>>\TT$
is denoted by $\TT_{>0}$. Thus $\io$ restricts to a bijection

(a) $\cb^2_{>0}@>\si>>\TT_{>0}$.
\nl
It follows that $\TT_{>0}$ is a cell of dimension $2\nu$.
We say that $\TT_{>0}$ is the totally positive part of $\TT$.

\head 2. A map from $G_{>0}$ to $\TT_{>0}$\endhead
\subhead 2.1\endsubhead
Let $g\in G_{>0}$. Now $g$ acts on the Lie algebra $\fg$ of $G$ by the
adjoint action. By \cite{L19, 9.1(a)} we have $\fg=\op_{a\in\RR_{>0}}\fg_a$
where $\fg_a$ is the $a$-eigenspace of $Ad(g):\fg@>>>\fg$. Let
$$\fg_{\ge1}=\op_{a\in\RR_{>0};a\ge1}\fg_a,$$
$$\fg_{\le1}=\op_{a\in\RR_{>0};a\le1}\fg_a.$$
In \cite{L21,\S5} it is shown that $\fg_{\ge1}$ is the Lie algebra
of a Borel subgroup $B$ of $G$, that  $\fg_{\le1}$ is the Lie algebra
of a Borel subgroup $B'$ of $G$ and that $\fg_1$ is the Lie algebra
of $B\cap B'$ which is a maximal torus of $G$; moreover we have
$g\in B\cap B'$. In \cite{L21,\S5} it is also shown that $B\in\cb_{>0}$.
The same argument applied to $g\i\in G_{<0}$ (for which the
$a$-eigenspace is equal to the $a\i$-eigenspace of $g$) shows that
$B'\in\cb_{<0}$. Thus $g\m(B,B')$ is a well defined map
$\p':G_{>0}@>>>\cb^2_{>0}$.

\subhead 2.2\endsubhead
Another approach to the map $\p'$ in 2.1 was pointed out to me by Xuhua
He, as an answer to my question.
Let $g\in G_{>0}$. By \cite{L94, 8.9} there is a unique $B\in\cb_{>0}$
such that $g\in B$. Similarly, by \cite{L94, 8.9} (applied to $G$ with
the new pinning as in 1.1 and to $g\i\in G_{<0}$) we see that there is
a unique $B'\in\cb_{<0}$ such that $g\i\in B'$. Then we have
$g\in B\cap B'$. The map $G_{>0}@>>>\cb^2_{>0}$, $g\m(B,B')$ is the same
as the map $\p'$ in 2.1 (this follows from results in \cite{L21,\S5}).

\subhead 2.3\endsubhead
Composing $\p':G_{>0}@>>>\cb^2_{>0}$ with the bijection
$\cb^2_{>0}@>\si>>\TT_{>0}$ in 2.5 we obtain a map $\p:G_{>0}@>>>\TT_{>0}$.

Using the arguments in 2.1 we see that $\p$
is the map which to any $g\in G_{>0}$ associates the
connected centralizer $Z_G^0(g)$ of $g$ in $G$. (From \cite{L94, 5.6} it
is known that $Z_G^0(g)\in\TT$.) In particular we must have
$Z_G^0(g)\in\TT_{>0}$ for any $g\in G_{>0}$.

\head 3. The case of $GL_2(\CC)$\endhead
\subhead 3.1\endsubhead
Assume that $G=GL_2(\CC)$ with $B^+$ (resp. $B^-$) being the group of
upper triangular (resp. lower triangular) matrices in $G$, and
$$x_1(a)=\left(\mat 1&a\\0&1\endmat\right),$$
$$y_1(a)=\left(\mat 1&0\\a&1\endmat\right)$$
for $a\in\CC$; here $I=\{1\}$.
In this case $G_{>0}$ consists of all $g\in G$ with all entries in $\RR_{>0}$
and with determinant in $\RR_{>0}$. We have $U^+_{>0}=\{x_1(a);a\in\RR_{>0}$,
$U^-_{>0}=\{y_1(a);a\in\RR_{>0}$.

Let $\ct\in\TT_{>0}$. Let $(B,B')\in\cb^2_{>0}$ be such that $\ct=B\cap B'$.
 We wish to describe the intersection
$\ct\cap G_{>0}=B\cap B'\cap G_{>0}$. We have $B=y_1(a)B^+y_1(-a)$,
$B'=y_1(-c)B^+y_1(c)$ where $a\in\RR_{>0}$, $c\in\RR_{>0}$.

By the proof of 1.2, $B\cap B'$ consists of all matrices of the form
$$y_1(a)x_1(-1/(a+c))\left(\mat t&0\\0&s\endmat\right)
x_1(1/(a+c))y_1(-a)$$
for various $t,s$ in $\CC^*$ that is of all matrices of the form
$$\left(\mat(tc+sa)/(a+c)&(t-s)/(a+c)  \\
ac(t-s)/(a+c)&(ta+sc)/(a+c)\endmat\right)$$
for various $t,s$ in $\CC^*$. The condition that such a matrix (with
$t,s$ in $\RR_{>0}$) is in $G_{>0}$ is that $t/s>1$.
We see that $\ct\cap G_{>0}$ is nonempty.

\head 4. The case of $GL_3(\CC)$\endhead
\subhead 4.1\endsubhead
Assume that $G=GL_3(\CC)$
with $B^+$ (resp. $B^-$) being the group of
upper triangular (resp. lower triangular) matrices in $G$, and
$$x_1(a)=\left(\mat 1&a&0\\0&1&0\\0&0&1\endmat\right),$$
$$x_2(a)=\left(\mat 1&0&0\\0&1&a\\0&0&1\endmat\right),$$
$$y_1(a)=\left(\mat 1&0&0\\a&1&0\\0&0&1\endmat\right),$$
$$y_2(a)=\left(\mat 1&0&0\\0&1&0\\0&a&1\endmat\right),$$
for $a\in\CC$; here $I=\{1,2\}$; $T$ consists of the diagonal matrices
$$\d(t,s,r)=\left(\mat t&0&0\\0&s&0\\0&0&r\endmat\right)$$
for various $(t,s,r)\in(\CC^*)^3$.

Now, $G_{>0}$ consists of all $g=(g_{ij})\in G$ such that
$g_{i,j}\in\RR_{>0}$ for all $i,j$ in $\{1,2,3\}$
and $g_{i,j}g_{i+1,j+1}-g_{i,j+1}g_{i+1,j}\in\RR_{>0}$
for $(i,j)\in\{(1,1),(1,2),(2,1),(2,2)\}$, $\det(g)>0$.
(The other minors of $g$ are then automatically in $\RR_{>0}$.)

For $$M=\left(\mat 1&0&0\\x&1&0\\z&y&1\endmat\right)\in U^-_{<0}$$
we have
$$\ti M=\left(\mat1&y/(xy-z)&1/z\\0&1&x/z\\0&0&1\endmat\right)\in U^+_{<0}.$$
Here $x,y$ are in $\RR_{<0}$, and $z,xy-z$ are in $\RR_{>0}$.
Let $\ct\in\TT_{>0}$.
By the proof of 1.2, there are well defined elements
$u=\left(\mat1&0&0\\a&1&0\\b&c&1\endmat\right)\in U^-_{>0}$,
$v=\left(\mat 1&0&0\\a'&1&0\\b'&c'&1\endmat\right)\in U^-_{>0}$
(here $a,b,c,a',b',c',ac-b,a'c'-b'$ are in $\RR_{>0}$) such that
$$\ct=\{S\d(t,s,r)S\i;(t,s,r)\in(\CC^*)^3\}$$ 
where
$$\align&S=u\wt{(u\i v\i)}=\left(\mat1&0&0\\a&1&0\\b&c&1\endmat\right)
\left(\mat1&0&0\\-a-a'&1&0\\A&-c-c'&1\endmat\right)\ti{}\\&
=\left(\mat1&0&0\\a&1&0\\b&c&1\endmat\right)
\left(\mat1&-(c+c')/C&1/A\\0&1&-(a+a')/A\\0&0&1\endmat\right)\\&
=\left(\mat1&-(c+c')/C&1/A\\a&(-ac+b+b')/C&-a'/A\\
b&(acc'-bc'+b'c)/C&(a'c'-b')/A\endmat\right)\endalign$$
and
$$\align&S\i=(\wt{(u\i v\i)})\i u\i\\&
=\left(\mat b'/C&c'/C&1/C\\ -(aa'c'-ab'+a'b)/A&
(a'c'-b-b')/A&(a+a')/A\\ac-b&-c&1\endmat\right).\endalign$$   
Here $A=ac+a'c+a'c'-b-b'\in\RR_{>0}$, $C=ac'+b+b'\in\RR_{>0}$.
For $(t,s,r)\in(\CC^*)^3$ we have $S\d(t,s,r)S\i=g=(g_{ij})$, where
$$g_{11}=b't/C+(c+c')(aa'c'-ab'+ab')s+(ac-b)r/A,$$
$$g_{12}=c't/C-(c+c')(a'c'-b-b')s-cr/A=c'(t-s)/C+c(s-r)/A,$$
$$g_{13}=t/C-(a+a')(c+c')s+r/A=(t-s)/C-(s-r)/A,$$
$$\align&g_{21}=ab't/C-(b+b'-ac)(aa'c'-ab'+a'b)s-a'(ac-b)r/A\\&
=ab'(t-s)/C+a'(ac-b)(s-r)/A,\endalign$$   
$$\align&g_{22}=ac't/C+(b+b'-ac)(a'c'-b-b')s+a'cr/A\\&
=ac'(t-s)/C+(a'c'-b'+ac-b)(b+b'+ac')s+a'cr/A,\endalign$$   
$$g_{23}=at/C+(b+b'-ac)(a+a')s-a'r/A=a(t-s)/C+a'(s-r)/A,$$
$$\align&g_{31}=bb't/C-(acc'-bc'+b'c)(aa'c'-ab'+a'b)s+(ac-b)(a'c'-b')r/A\\&
=bb'(t-s)/C-(ac-b)(a'c'-b')(s-r)/A,\endalign$$   
$$\align&g_{32}=bc't/A+(acc'-bc'+b'c)(a'c'-b-b')s-c(a'c'-b')r/A\\&
=bc'(t-s)/A+c(a'c'-b')(s-r)/A,\endalign$$   
$$g_{33}=bt/C+(acc'-bc'+b'c)(a+a')s+(a'c'-b')r/A.$$
We have $g\i=S\i\d(t\i,s\i,r\i)S=(g'_{ij})$, where
$$g'_{11}=b't\i/C+ac's\i/C+br\i/C,$$
$$\align&g'_{13}=b't\i/(AC)-a'cs\i/(AC)+(a'c'-b')r\i/(AC)\\&
=(a'c'-b')(r\i-s\i)/(AC)-b'(s\i-t\i)/(AC),\endalign$$   
$$g'_{31}=(ac-b)t\i-acs\i+br\i=b(r\i-s\i)-(ac-b)(s\i-t\i)$$
$$g'_{33}=(ac-b)t\i/A+a'cs\i/A+(a'c'-b')r\i/A$$
(we do not write the other $g'_{ij}$).

\subhead 4.2\endsubhead
The condition that $g\in G_{>0}$ is that $g_{ij}\in\RR_{>0}$ for all $i,j$ and
that

$g'_{11},g'_{13},g'_{31},g'_{33}$
\nl
are in $\RR_{>0}$. (We must also have
$(t,s,r)\in\RR_{>0}^3$, hence $\det(g_{ij})=tsr>0$.)
From $g_{13}\in\RR_{>0}$ we see that if $s-r\ge0$ then $t-s>0$.
(Hence if $t-s\le0$ then $s-r<0$.)
From $g'_{31}\in\RR_{>0}$ we see that if $t-s\ge0$ then $s-r>0$.
(Hence if $s-r\le0$ then $t-s<0$.) We see that we have
either $t>s>r$ or $t<s<r$. But $t<s<r$ would contradict $g_{12}\in\RR_{>0}$.
Thus
we must have $t>s>r$. It follows that the condition that

$g_{11},g_{12},g_{21},g_{22},g_{23},g_{32},g_{33}, g'_{11},g'_{33}$
\nl
are in $\RR_{>0}$ is
automatically satisfied. The remaining conditions can be written as

$(t-s)/(s-r)>C/A, (t-s)/(s-r)>(ac-b)(a'c'-b')C/(bb'A)$,

$(r\i-s\i)/(s\i-t\i)>b'/(a'c'-b'),(r\i-s\i)/(s\i-t\i)>(ac-b)/b$

or equivalently

(a) $(t/s-1)/(1-r/s)>C/A, (t/s-1)/(1-r/s)>(ac-b)(a'c'-b')C/(bb'A)$,

(b) $(s/r-1)/(1-s/t)>b'/(a'c'-b'),(s/r-1)/(1-s/t)>(ac-b)/b$.
\nl
We have $1/(1-r/s)>1,1/(1-s/t)>1$ so that $(t/s-1)/(1-r/s)>t/s-1$,
$(s/r-1)/(1-s/t)>s/r-1$. Thus if $t/s$ and $s/r$ are sufficiently large,
then (a),(b) hold.

We see that the following holds.

\proclaim{Proposition 4.3}$\ct\cap G_{>0}$ is equal to the set of all
$S\d(t,s,r)S\i$ with $(t,s,r)\in\RR_{>0}^3$ satisfying  $t>s>r$ and
4.2(a),(b). This set contains all $S\d(t,s,r)S\i$ with
$(t,s,r)\in\RR_{>0}$ such that $t/s$ and $s/r$ are sufficiently large.
In particular, it is nonempty.
\endproclaim

\subhead 4.4\endsubhead
Let $u,a,b,c$ be as in 4.1. Then $B=uB^+u\i\in\cb_{>0}$.
Let $\ce_B$ be the set of all $(t,s,r)\in \RR_{>0}^3$ which can appear
as eigenvalues of some element in $B\cap G_{>0}$ (we can assume that
$t>s>r$). From 4.3 it follows that $\ce_B$ consists of all $t>s>r>0$ such
that 4.2(a),(b) hold for some $a',b',c'$ in $\RR_{>0}$ such that
$a'c'-b'>0$.
In particular any $(t,s,r)\in\ce_B$ must satisfy
$(s/r-1)/(1-s/t)>(ac-b)/b$. We see that $\ce_B$ can be a proper
subset of $\{(t,s,r)\in\RR_{>0}^3;t>s>r\}$.

\head 5. A conjecture\endhead
\subhead 5.1\endsubhead
We return to the setup in 0.1.
For $i\in I$ we define $\c_i:T@>>>\CC^*$ by
$tx_i(a)t\i=x_i(\c_i(t)a)$ for all $t\in T,a\in\CC$.
For any $p\in\RR_{>0}$ let
$$T^p_{>0}=\{t\in T_{>0};\c_i(t)>p\text{ for all }i\in I\}.$$

Let $\ct\in\TT_{>0}$. We have $\ct=STS\i$ where
$S=u\wt{(u\i v\i)}$ for well defined elements $u,v$ in $U^-_{>0}$.
We conjecture that 

(a) $\ct\cap G_{>0}\sub ST^1_{>0}S\i$ and

(b) $ST^p_{>0}S\i\sub\ct\cap G_{>0}$ for some $p\ge1$.
\nl
In particular this would imply that $\ct\cap G_{>0}\ne\emp$
(so that the map $\p$ in 2.3 is surjective).

I expect that (b) can be proved by a variant of the argument in
\cite{L94, 7.1}.

This conjecture holds when $G$ is $GL_2$ or $GL_3$ (see \S3,\S4).

\enddocument